\documentclass[12pt]{amsart}

\usepackage{amsfonts,amscd,amsmath,amssymb,amsthm,latexsym,xypic}
\newcommand{\pf}{\vspace{\baselineskip}\noindent{\bf Proof.}$\;\;$}

\newcommand{\ZZ}{{\mathbb Z}}

\newcommand{\RR}{{\mathbb R}}

\newcommand{\n}{\noindent}
\newcommand{\NN}{{\bf N}}

\newtheorem{df}{Definition}[section]
\newtheorem{teo}[df]{Theorem}
\newtheorem{prop}[df]{Proposition}

\begin{document}
\author{Federica Galluzzi}
\address{Dipartimento di Matematica, Universit\`a di Torino, Via Carlo Alberto n.10, 10123 Torino, ITALY}
\email{federica.galluzzi@unito.it}

\thanks{The first author is supported by Progetto di Ricerca Nazionale COFIN 2004 "Geometria sulle Variet\`a
Algebriche"}

\author{Giuseppe Lombardo}
\address{Dipartimento di Matematica, Universit\`a di Torino, Via Carlo Alberto n.10, 10123 Torino, ITALY}
\email{lombardo@dm.unito.it}
\title{On the automorphisms of some $K3$ surface double cover of the plane}
\maketitle

\begin{abstract}
In this paper we study the automorphisms group of some $K3$ surfaces
which are  double covers of the projective plane ramified over a
smooth sextic plane curve. More precisely, we study the case of a
$K3$ surface of Picard rank two such that there is a rational curve of degree $d$ which is tangent to the sextic in $d$ points.
\end{abstract}

\medskip
\noindent \centerline{{\bf Introduction}}

\n 
$K3$ surfaces which are double covers of the plane ramified over a plane sextic are classical objects. In this paper we determine the automorphisms group of some of these surfaces. More precisely, we restrict to the case of Picard rank two. Moreover we suppose
that there is a rational curve of degree $d$ which is tangent to
the sextic in $d$ points. The result is that $K3$ surfaces of this type with $d$ odd have automorphisms group isomorphic to $\ZZ _2.$

\n
We first analize the case $d=3 \,$ in full details to better explain the results and the methodology we have used. Then we study the 
case with $d$ odd for which the situation is not substantially different from the previous one.

\n The automorphisms of a $K3$ surface are given by the Hodge
isometries that preserve the K\"ahler cone (see \cite{BPV},
VIII.11). Thus, the strategy to study the automorphism of a $K3$
surface $X$ is to determine its K\"ahler cone and the Hodge isometries of
$H^2(X,\ZZ)$ which preserve it. To do this, we first determine the isometries of the
N\'eron-Severi lattice  $NS(X)$ which preserve the K\"ahler cone and then we study the gluing conditions on the transcendental lattice $T(X),\,$ since $H^2(X,\ZZ)=NS(X)\oplus T(X).\,$ 

\n
In the preliminary Section \ref{prel} we introduce some basic material on lattices
and $K3$ surfaces. In section \ref{xd} we give explicitly the N\'eron-Severi lattice of the $K3$  we are interested in. In section \ref{d=3} we start our work on the case $d=3$ following the strategy explained before. We determine the K\"ahler cone in Prop.\ref{cone} and the automorhisms group of the N\'eron-Severi lattice in Prop.\ref{l3}  using some basic facts on the generalized Pell equation $a^2 - 13 b ^2 = -4.\,$ Finally, we derive in Theorem\ref{auto3} that the automorhisms group of the surface is isomorphic to $\ZZ_2$ analyzing the gluing condition on the transcendental lattice. Using similar methods we obtain the same result 
for $d$ odd in Theorem\ref{teodd}.

\medskip
\n
\section{Preliminaries}\label{prel}\subsection{Lattices}
A {\em lattice} is a free $\ZZ$-module $L$ of finite rank  with a $\ZZ$-valued symmetric bilinear form $<,>.\,$ A lattice  is called {\em even} if the quadratic form associated to the bilinear form has only even values, {\em odd} otherwise. The {\em discriminant} $d(L)$ is the determinant of the matrix of the bilinear form. A lattice is called {\em non-degenerate} if the discriminant is non-zero and {\em unimodular} if the discriminant is $\pm 1 .\,$ If the lattice $L$ is non-degenerate, the pair $(s_+, s_-),\,$ where $s_{\pm}$ denotes the multiplicity of the eigenvalue $\pm 1$ for the quadratic form associated to $L \otimes \RR ,\,$ is called {\em signature} of $L .\,$ Finally, we call $\, s_+ +s_- \,$ the {\em rank} of $L.$

\n
Given a lattice $(L,<,>) \,$ we can construct the lattice $(L (m) <,>_m),\,$ that is the $\ZZ$-module  $L $ with form $<x,y>_m=m<x,y>.$

\n
An {\em isometry} of lattices is an isomorphism preserving the bilinear form.
Given a sublattice $L \hookrightarrow L ^{\prime},\,$ the embedding is {\em primitive} if $\displaystyle {L ^{\prime} \over L }$
is free. Two even lattices $S,\, T$ are \textit{orthogonal} if there exist an even unimodular lattice $L$ and a primitive embedding $S \hookrightarrow L$ for which $(S)_L^{\bot} \cong T.$

\subsection{$K3$ surfaces}
\n A $K3$ surface is a compact K\"ahler surface with trivial canonical bundle and such that its first Betti number is equal to zero. Let $U$ be the lattice of rank two with quadratic form given by the matrix $ \begin{pmatrix}0&1 \cr 1&0 \end{pmatrix}\,$ and 
let $E_8$ be the lattice of rank eight whose quadratic form is the  Cartan matrix of the root system of $E_8$. It is an even, unimodular and positive definite lattice.

\n
It is well known that $H^2(X,\ZZ)$
is an even lattice of rank $22$ and signature $(3,19)$ isomorphic
to the lattice
$$
\Lambda = U \oplus U \oplus U \oplus E_8 (-1) \oplus E_8 (-1),
$$
that we will call, from now on, the $K3-$\textit{lattice}.
\n Denote with $NS(X) \cong H^2(X,\ZZ) \cap H^{1,1}(X)$
the N\'eron-Severi lattice of $X$ and with $T(X)$ the orthogonal
complement of $NS(X)$ in $H^2(X,\ZZ).$ The \textit{Picard rank
of $X$},  $\rho (X),$ is the rank of $NS(X).$ The Hodge Index Theorem implies that
$NS(X)$ has signature $(1,\rho(X)-1)$ and that $T(X)$ has signature $(2,20-\rho(X)).\,$

\smallskip
\noindent
We will use the following result:

\n
\begin{teo}\cite[Thm. 1.14.4]{N}\cite[2.9]{M}\label{unique}
If $\rho(X) \leq 10 ,\,$then every even lattice $S$ of signature $(1,\rho-1 )$ occurs as the N\'eron-Severi group of some $K3$ surface and the primitive embedding $S \hookrightarrow \Lambda$ is unique.
\end{teo}

\medskip
\n
Denote with $\Delta$ the set of
the classes of the $(-2)$-curves in $NS(X)$ and with
$\mathcal C \subset NS(X)\otimes \RR$ the connected
component of the set of elements $x\in NS(X)\otimes \RR$ with $x ^
2>0 $ which contains an ample divisor.
The \textit{K\"ahler cone} is the convex subcone of $\mathcal C$ defined as
$$
\mathcal C^+=\{y \in \mathcal C\, :\, (y,d)>0\; \mbox{\, for all}\; d \in NS(X),\,d \ \mbox{effective} \}.
$$

\n
We will also use the following
\begin{prop}\label{kcone} \cite[VIII 3.8.]{BPV}
The K\"ahler cone is given by
$$
\mathcal C^+= \{  w  \in \mathcal C : \, w N
>0 , \, \mbox{for all} \ N \in \Delta \}.
$$
\end{prop}

\section{The surfaces $X_d$}\label{xd}

\n It is well known that a surface which is a double cover of the
projective plane ramified over a smooth sextic plane curve is a
$K3$ surface.  We restrict to the case when the N\'eron-Severi lattice has rank two
and such that there is a rational curve of degree $d$ which is tangent to the sextic in $d$ points. Let $X_d$ be such a surface. We suppose that $X_d$ is
general.
The N\'eron-Severi lattice of $X_d$ has  quadratic form given by
$$
Q_d:=\begin{pmatrix}
2 & d \\
d & -2
\end{pmatrix}.
$$

\n We denote this lattice $L_d.\,$ It is has segnature $(1,1).\,$

\n
Our aim is to compute the automorphisms group of
$X_d\,$ We start with the case $d=3.$

\section{Case $d=3$}\label{d=3}

\n  We want to study now the automorphisms group of a $K3$ surface $X_3$ 
of rank two which is a double cover of the plane ramified over a smooth sextic which has a rational tritangent cubic.
Such a surface has a N\'eron-Severi lattice $L_3$ given by the matrix
$$
Q_3:=\begin{pmatrix}
2 & 3 \\
3 & -2
\end{pmatrix}.
$$

\smallskip
\n
\subsection{The K\"ahler cone}

\n
Denote with $\mathcal C^+_3$ the K\"ahler cone of $X_3.\,$ We have the following

\smallskip
\n
\begin{prop}\label{cone}
Let $X_3$ be a surface with N\'eron-Severi lattice isomorphic to $L_3.\,$ Then there is an isomorphism of lattices $NS(X_3)\otimes \RR \cong L_3 \otimes \RR \cong \RR ^2$ such that:
$$
\mathcal C^+_3=\left \{(x,y) \in \RR ^2 \, : \,3x-2y >0 \right \} \cap \left \{ (x,y) \in \RR ^2 \, : \, 3x+11y >0\right \}.
$$
\end{prop}

\pf  We have first to determine the classes of the $(-2)-$curves on
$X_3\,$ that is the set $\Delta \subset NS(X_3)$ :
$$
\Delta = \{ D \in NS(X_3)\, : \, D > 0 , D^2=-2, D \ \mbox{irreducible}
 \}.
$$

\n The condition $D^2 = -2$ means that we have to determine the
integer solutions of the equation
\begin{equation}\label{eq1}
x^2+3xy-y^2=-1.
\end{equation}

\n
We write $x^2+3xy-y^2=(x- \alpha y)(x- \bar
{\alpha} y),\,$ with $\alpha = \frac{-3+\sqrt{13}}{2},\,$ and
$\bar{\alpha} = \frac{-3-\sqrt{13}}{2}=-3-\alpha .\,$ 
Thus $\Delta$ corresponds to the set
$$
\{u \in \ZZ[\alpha]\,: \,u \bar u =-1\}.
$$

\n
It is known that the invertible elements in
$\ZZ [\alpha]$ are $\ZZ [\alpha]^*=<\eta>$ where $\eta =\frac {3 + \sqrt{13}}{2}= \alpha +3\,$ and $\eta \bar{\eta}=-1.\,$ Thus, solutions of (\ref{eq1}) are given by the odd
powers of  $\eta \,$ and $\bar{\eta}.\,$ The element $\eta$ verifies $\eta ^3= 11\eta + \bar {\eta},\,$  so by induction we
obtain that  $\,\eta ^{2k+1} = a \eta +
b \bar{\eta},\,$ where $a,b \in \ZZ _{>0}.\,$ This shows that to $\eta$ and $\bar{\eta}$ correspond the unique two
irreducible $(-2)$-curves $D_{\eta},\, D_{\bar{\eta}}$ in $\Delta .\,$

\n
The element $\eta$ represents the solution $(0,1)$ of the equation (\ref{eq1}) and $\bar{\eta}=3 -\eta$ represents $(3,-1).\,$
Now, we
can determine $\mathcal C^+$ that is,  following Prop.\ref{kcone}, the set
$$
\mathcal C^+= \{  w  \in \mathcal C : \, Q_3(w, D)
>0 , \, \mbox{for all} \ D \in \Delta \}.
$$

\n
This means that we are looking for the elements $w=(x,y)$ such that $Q_3(w, \eta)>0$ and $Q_3(w,\bar{\eta} )>0,\,$
thus we obtain the statement.
\qed

\medskip
\n
\subsection{The automorphisms group}

\n  Denote with $T_3$ the transcendental lattice of $X_3.\,$ As $H^2(X_3,\ZZ)\cong L_3 \oplus T_3,\,$ we start studying the isometries of $L_3,\,$ then we'll identify the ones preserving the ample cone and finally  we'll analyze the gluing conditions on $T_3.\,$

\smallskip
\noindent
First,
we note that there exists a primitive embedding  $L_3 \overset{i}{\hookrightarrow} U
\oplus U ,\,$ but it is not
unique. So one can have another embedding $L_3 \overset{j}{\hookrightarrow} U \oplus U$ such that $i(L_3)^{\bot}\ncong j(L_3)^{\bot}\,$ in $U \oplus U.\,$ In any case, by Theorem \ref{unique} one has
$$
i(L_3)^{\bot} \oplus U \oplus E_8(-1) \oplus E_8(-1) \cong j(L_3)^{\bot} \oplus U \oplus E_8(-1)\oplus E_8(-1).
$$

\n
Moreover, $L_3$ is orthogonal to $L_3(-1)$ in $U \oplus U.$
For this reason we can choose $T_3= L_3.\,$

\newpage
\n
\begin{prop}\label{l3}  The automorphisms group $Aut(L_3)$ is isomorphic to
the group
$
\ZZ _2  *  \ZZ _2 .
$
\end{prop}

\pf The group of isometries of $L_3$ are given by
$$
O(L_3)= \left \{ M \in GL_2(\ZZ)\, :\, ^t\!M \, Q_3\,  M = Q_3 \right \}
$$

\n By direct computations one obtains matrices of the following
form

\smallskip
$$
P^{\pm}_{(a,b)}:=\; \begin{pmatrix}
\displaystyle{\frac{11b \mp 3a}{2} } &\displaystyle{\frac{-3b \pm a}{2} } \\
\displaystyle{\frac{-3b \pm a}{2}}&b\\
\end{pmatrix} ,\;\;\;\;
Q^{\pm}_{(a,b)} := \begin{pmatrix}
- b &\displaystyle{\frac{-3b \mp a}{2} } \\
\displaystyle{\frac{-3b \pm a}{2}}&b\\
\end{pmatrix}.
$$

\smallskip
\n where the $(a,b)$ are solutions of the generalized Pell equation
\begin{equation}\label{pell4}
a^2 - 13 b ^2 = -4
\end{equation}

\n
A standard result on Pell equations and on fundamental units, see for example \cite{Ba} and \cite{Co}, says that the solutions  are $(\pm a_n,\pm b_n)$ with  $\displaystyle{\frac{a_n + \sqrt{13} b_n}{2}}=\left (\frac{a_0+\sqrt{13}\ b_0}{2}\right)^{2n+1},\,n \in \NN$ and $(a_0,b_0)$ is the pair of smallest positive integers satisfying the equation.
In our
case the pair of smallest positive integers that satisfy the Pell
equation (\ref{pell4}) is $(a_0,b_0)=(3,1).\,$ Notice that $\displaystyle{\frac{a_0+\sqrt{13}\ b_0}{2}=\eta}$ and then we can obtain solutions
$(a_n,b_n)$ by recurrence multiplying by $\eta ^2.\,$ By direct computations:
$$
\begin{cases}
\displaystyle{a_{n+1}=\frac{11a_n + 39 b_n}{2}}\\
\displaystyle{b_{n+1}=\frac{3a_n + 11 b_n}{2}}\; \;\;.\\
\end{cases}
$$

\smallskip
\n
Moreover, if $(a_n,b_n)$ gives rise to the matrices  $ P^{\pm}_{(a_{n},b_{n})},\,Q^{\pm}_{(a_{n},b_{n})},\,$ then
the couples $(a_n,-b_n),\,(-a_n,b_n),\,(-a_n,-b_n)$ give rice to the matrices
$$
(- P^{\mp}_{(a_n,b_n)},\,-Q^{\pm}_{(a_n ,b_n)}),\,(P^{\mp}_{(a_n,b_n)},\,-Q^{\pm}_{(a_n,b_n)}),\,(- P^{\pm}_{(a_n,b_n)},\,-Q^{\mp}_{(a_{n},b_{n})})
$$

\n
respectively.
We write $P^{\pm}_{n}:=P^{\pm}_{(a_{n},b_{n})}$ and $Q^{\pm}_{n}:=Q^{\pm}_{(a_{n},b_{n})}.\,$
For $(3,1)$ one obtains the matrices
$$
P^+_{0}= I\,,\;\;\;\;P^-_{0}=\;\begin{pmatrix}
10& -3 \\
-3&1\\
\end{pmatrix} ,\;\;\;\;
Q^+_{0}=\begin{pmatrix}
-1&0 \\
-3&1\\
\end{pmatrix}
\;\;\;\;
Q^-_{0}=\begin{pmatrix}
-1&-3 \\
0&1\\
\end{pmatrix}.
$$

\n The matrices $Q^+_{0},\,Q^-_{0}$ are non commuting involutions and $P^-_{0}=Q^-_{0}Q^+_{0}.\,$ The matrices $P^{\pm}_{n+1},\,Q^{\pm}_{n+1}$ are obtained by multiplication
$$
P^+_{n}=(P_1^+)^n=(P_0^-)^{-n},\;\;\; P^-_{n}= (P_0^-)^{n+1}
,\;\;\;\;
Q^+_{n+1}=P_0^- Q^+_{n},\;\;\; Q^-_{n+1}= Q^-_{n} P_0^-
$$

\n
and $(P_0^-)^{-1}=P_1^+.\,$
Set $p$ and $q$ for the automorphism of $L_3$ corresponding to
$ Q^+_{0}$ and $Q^-_{0}$ respectively. Thus we have showed that the group $O(L_3)$ can
be described as $\langle p \rangle * \langle q \rangle .\,$  \qed

\bigskip
\n
\begin{teo}\label{auto3}
The automorphisms group of $X_3$ is isomorphic to $\ZZ _2.$
\end{teo}

\pf We are looking for Hodge isometries of $H^2(X_3,\ZZ)$ which
preserve the ample cone. From the generality of
$X_3$  we may assume that the only Hodge isometries of
$T_3$ are $\pm I.\,$ Thus, we have first to identify the elements
in $Aut(L_3)$ which preserve the K\"ahler cone and then we impose
a gluing condition on $T_3,\,$ since the isometries we are looking
for have to induce $\pm I$ on $T_3.$ Note first that $-I \in Aut(L_3)$
can not preserve the ample cone.
We have from Prop.\ref{cone} that the K\"ahler cone is isomorphic to the chamber delimited by $H_{D_{\eta}}\cong \RR^+\!\langle \mathbf{v} \rangle$ and $\,H_{D_{\bar{\eta}}}\cong \RR^+\!\langle \mathbf{w} \rangle$ where $\mathbf{v}=(2,3)$ and $\mathbf{w}=(11,-3).\,$

\noindent
An easy computation shows that
the elements in $Aut(L_3)$ having this property are the ones
generated by $-q \,$  which forms a $\ZZ _2.\,$ A direct
computations gives that $-q$ satisfy the gluing condition on
$T_3.\,$ \qed

\section{Case  $d$ odd.}\label{d=odd}

\n In this case we have a $K3$ surface with N\'eron-Severi lattice
$L_d$ of rank two given by the matrix
$$
Q_d:=\begin{pmatrix}
2 & d \\
d & -2
\end{pmatrix}.
$$

\smallskip
\n Set  $X_d$ for the $K3$ surface having $NS(X_d)\cong L_d.\,$
Such a $K3$ is a double cover of the plane ramified over a smooth sextic 
tangent to a rational curve of degree d. 
Following the same strategy adopted for the case $d=3,\,$ we obtain

\begin{teo}\label{teodd}
If $d$ is odd the automorphisms group of $X_d$ is isomorphic to $\ZZ _2.$
\end{teo}

\pf We compute
$$
O(L_d)= \left \{ M \in GL_2(\ZZ)\, :\, ^t\!M \, Q_d\,  M = Q_d
\right \}
$$

\n and we obtain matrices of the following form

\smallskip
$$
R^{\pm}_{(a,b)}:=\; \begin{pmatrix}
\displaystyle{\frac{(2+d^2)b \mp da}{2} } &\displaystyle{\frac{-db \pm a}{2} } \\
\displaystyle{\frac{-db \pm a}{2}}&b\\
\end{pmatrix} ,\;\;\;\;
S^{\pm}_{(a,b)} := \begin{pmatrix}
- b &\displaystyle{\frac{-db \pm a}{2} } \\
\displaystyle{\frac{-db \mp a}{2}}&b\\
\end{pmatrix}.
$$

\smallskip
\n where the $(a,b)$ are solutions of the Pell equation
\begin{equation}\label{genpell}
a^2 - (d^2+4) b ^2 = -4.
\end{equation}

\n When $d$ is odd, by theory on Pell equation the situation is analogous to the one of Prop.\ref{l3} that is, all solutions can be generated from the minimal positive solution.

\n This means that $\ZZ [\sqrt{d^2+4}]^*$ is generated by $\displaystyle{\eta =\frac{d}{2}+ \frac{\sqrt{d^2+4}}{2}}$
and that if $(a_n,b_n)$ is a solution of (\ref{genpell}),  then
$(a_{n+1},b_{n+1})$ is obtained by multiplying by $\eta ^2.\,$ By
direct computations, the solutions are obtained by recurrence
$$
\begin{cases}
\displaystyle{a_{n+1}=\frac{a_nd^2 +2 a_n+b_n d^3 +4 b_nd}{2}}\\
\displaystyle{b_{n+1}=\frac{a_n d +  b_n d^2+ 2 b_n}{2}}\; \;\;.\\
\end{cases}
$$

\smallskip
\n The pair of smallest positive integers that satisfy the Pell
equation (\ref{pell4}) is $(a_0,b_0)=(d,1).\,$

\n
We write $R^{\pm}_{n}:=R^{\pm}_{(a_{n},b_{n})}$ and
$S^{\pm}_{n}:=S^{\pm}_{(a_{n},b_{n})}.\,$ For $(d,1)$ one obtains
the matrices
$$
R^+_{0}= I\,,\;\;\;\;R^-_{0}=\;\begin{pmatrix}
1+d^2& -d\\
-d&1\\
\end{pmatrix} ,\;\;\;\;
S^+_{0}=\begin{pmatrix}
-1&0 \\
-d&1\\
\end{pmatrix}
\;\;\;\; S^-_{0}=\begin{pmatrix}
-1&-d \\
0&1\\
\end{pmatrix}
$$

\n  The matrices
$S^+_{0},\,S^-_{0}$ are non commuting involutions and
$R^-_{0}=S^-_{0}S^+_{0}.\,$ 
The matrices
$R^{\pm}_{n+1},\,S^{\pm}_{n+1}$ are obtained by multiplication
$$
R^+_{n}=(R^+_1)^n=(R_0^-)^{-n} ,\;\;\;\; S^+_{n+1}=R_0^-S_n^+,\;\;\;S^-_{n+1}=S_n^-R_0^-
$$

\n
and $(R_0^-)^{-1}=R_1^+.\,$

\n Set $r$ and $s$ for the automorphism of $L_3$ corresponding to
$ S^+_{0}$ and $S^-_{0}$ respectively. The group $O(L_d)$ can be
described as $\langle r \rangle * \langle s \rangle .\,$  
We obtain that the K\"ahler cone is isomorphic to
$$
\mathcal C^+_d=\left \{(x,y) \in \RR ^2 \, : \,dx-2y >0 \right \} \cap \left \{ (x,y) \in \RR ^2 \, : \, dx+(d^2+2)y >0\right \}
$$

\n
and, as before, the only automorphism of the N\'eron-Severi lattice which preserves the cone is $-s\,$ and it satisfies the gluing condition on $T_d.$

\qed

\centerline{\textsc{Acknowledgments}}

\noindent
We would like to thank Prof. Bert van Geemen for the very useful discussions and valuable suggestions.

\end{document}